\documentstyle[draft,amscd,syntonly,amssymb]{amsart}

\theoremstyle{plain}
\newtheorem{Thm}{Theorem}
\newtheorem{Prop}[Thm]{Proposition}
\newtheorem{Cor}[Thm]{Corollary}
\newtheorem{Lem}[Thm]{Lemma}

 \theoremstyle{definition}
\theoremstyle{remark}

\numberwithin{equation}{section}

\begin{document} 
 \title{Symplectic action around loops in $\text{Ham}(M)$}

 \author{ ANDR\'{E}S   VI\~{N}A}
\address{Departamento de F\'{i}sica. Universidad de Oviedo.   Avda Calvo
 Sotelo.     33007 Oviedo. Spain. } 
 \email{avinae@@correo.uniovi.es}
\thanks{The author was partially supported by Universidad de Oviedo, grant MB-02-514}

 \maketitle
\begin{abstract}

Let $\text{Ham(M)}$ be the group of Hamiltonian
 symplectomorphisms of  a quantizable, compact, symplectic manifold $(M,\omega)$. 
We prove the existence of an action integral around loops in $\text{Ham(M)}$,  
 and determine the 
value of this action integral  on particular loops
when the manifold is a   coadjoint orbit.  
  \end{abstract}

\section {Introduction} \label{S:intro}

Given a compact, symplectic manifold $(M,\omega)$, the group of Hamiltonian symplectomorphisms
\cite{Mc-S} \cite{lP01}
of $(M,\omega)$
is denoted by $\text{Ham}(M)$.    If $\psi=\{\psi_t\}_{t\in[0,1]}$ is 
 a loop in $\text{Ham}(M)$ at $\text{id}$ and $x$ is a point of $M$,
then the closed curve in $M$ $\{\psi_t(x) \}_t$ is nullhomologous \cite{Mc-S} \cite{L-M-P}. 
If $S$ 
is a $2$-chain in $M$
whose boundary is this curve, one can consider 
\begin{equation}\label{kappa}
\kappa_x (\psi):=\text{exp}\,  \Big(2\pi i\int_S \omega-2\pi i\int_0^1 f_t(\psi_t(x)) dt \Big),
\end{equation}
where $f_t$ is the normalized time-dependent Hamiltonian associated to $\psi$.
Throughout this article we will assume that 
 $(M,\omega)$ is {\em quantizable} \cite{nW92}, that is, $\omega$ defines an integral cohomology class in $M$,
then the right hand side of (\ref{kappa}) is independent of the $2$-chain $S$.
 $\kappa_x(\psi)$
is in fact the $U(1)$-valued {\em action integral} around 
the curve $\{\psi_t(x) \}_t$ \cite{Mc-S} \cite{aW89}.
In \cite{aV01} we proved that $\kappa_x(\psi)$ is independent of the point $x$; the 
proof is based in an analysis
of some properties of the prequantization representation \cite{jS80}. Here we give a new
proof of this property in the context 
of the gauge transformations of a prequantum bundle  \cite{nW92}.
The idea of the proof is the following:
Since $M$ is quantizable, there is a prequantum bundle  over $M$, that is,
a Hermitian line bundle $\pi:L\rightarrow M$ with a connection  $D$, such that the 
curvature of $D$ is $-2\pi i\omega$.  The  time-dependent Hamiltonian $f_t$
determines the corresponding Hamiltonian vector fields $X_t$, and by $X^{\sharp}_t$
is denoted the horizontal lift of 
$X_t$. From  $f_t$ one can define a vertical vector field  $W_{f_t}$
whose value at $p$ is determined by the curve
 $\{p\cdot\text{exp}(2\pi itf(\pi(p)) \}_t$. Then $\{Z_t:=X^{\sharp}_t-W_{f_t}\}$ 
generates a family $F_t$ of 
diffeomorphisms of $L$ which preserve the connection. 
That is,  loops in $\text{Ham}(M)$ lift to preserving-connection 
isotopies of a prequantum bundle.
We prove that
 $F_1$ is the gauge transformation defined by the map
$x\in M\mapsto \kappa_x(\psi)\in U(1)$. It follows from this result that the map
$\kappa_{-}(\psi)$ is constant. From the independence of $\kappa_x(\psi)$ of $x$
one deduces that $\kappa(\psi)$ depends only on the homotopy class of $\psi$. So
we obtain a representation $\kappa$ of the group $\pi_1(\text{Ham}(M))$, which can be called
the action integral representation. The existence of this representation
has been proved in other contexts by Seidel \cite{pS97} and Schwarz \cite{mS00}.

The isotopy $F_t$ of $L$, lifting of the loop  $\psi$ in $\text{Ham}(M)$,
 allows us to assign to each section $\tau$ of $L$
a family $\tau_t$ of sections defined by
$\tau_t=F_t\circ\tau\circ\psi_t^{-1}$.
One can consider  the correspondence
$\tau\rightarrow\tau_t$ as   a ``transport"  of sections
of $L$ along $\psi$.  On the other hand,
 in the prequantization   of the manifold $M$ one constructs for  each 
 Hamiltonian vector field $X$  an operator ${\cal P}_X$
(the corresponding prequantization operator \cite{jS80}), which acts on the space $\Gamma(L)$
of sections of $L$.
If $X_t$ is the family of vector fields that determines the isotopy $\psi$, we prove that 
$$\frac{d\tau_t}{dt}={\cal P}_{X_t}(\tau_t),\;\; \; \tau_0=\tau$$
is the differential equation of the transport along $\psi$.  
In this context $\kappa(\psi)$ is the holonomy of this transport.

If $M$ is  a Hamiltonian $G$-space \cite{vG94} one can consider
 loops $\varphi$ in $\text{Ham}(M)$
defined by the $G$-action.
If a  family $g_t$ of elements in $G$ determines a loop $\varphi$ in $\text{Ham}(M)$,
the independence   of $\kappa_x(\varphi)$ from $x$ allows us to express
$\kappa(\varphi)$   in terms of the vertices of moment map. 
This fact, when the group $G$ is a torus $T$, implies the existence of
 a moment map for the $T$-action
whose vertices are integer lattice points of ${\frak t}^*$. In a Hamiltonian
$G$-space to each  $A\in{\frak g}$ corresponds a vector field $X_A$ in $M$.
The map 
${\cal P}$, which assigns to  each $A\in{\frak g}$ the operator ${\cal P}_{X_A}$,
is in fact a Lie algebra representation.  
The $G$-action is called 
{\em pre-quantizable} \cite{vG820}, if there is an action $\rho$ of $G$ on $L$ that induces
${\cal P}$. Guillemin and Sternberg proved the above  property about the vertices of a
   moment map under the
additional hypothesis that the $G$-action is
 pre-quantizable \cite[$\S$ 8, Corollary 1]{vG820}.

In the particular case, when the $G$-action is pre-quantizable,
we will express $\kappa(\varphi)$ in terms of
 the character of the restriction of $\rho$
to invariant  finite-dimensional subspaces of $\Gamma(L)$. More precisely,  
the action $\rho$ induces a representation $\upsilon$ of $G$ on the space 
$\Gamma(L)$. On the other hand, if
  $G$ is a compact group there are $G$-invariant almost complex structures on $M$.
Let $I$ be such an almost complex structure and ${\cal Q}_I$ the subspace of $\Gamma(L)$
consisting of the $I$-polarized sections. Using that $\kappa(\varphi)$
is the holonomy of the transport along $\varphi$,   $\kappa(\varphi)$ can be 
expressed in terms of  the character of $\upsilon$ restricted 
to ${\cal Q}_I$ (Theorem \ref{Thchi1}).

  The coadjoint orbits of a compact Lie group $G$ are particular cases of 
  Hamiltonian $G$-spaces. If $\eta\in{\frak g}^*$ is integral, i.e. there is a character $\Lambda $ of the stabilizer
subgroup $G_{\eta}$ of $\eta$ whose derivative is $2\pi i\eta$, then the action of $G$ on
the orbit of $\eta$, ${\cal O}_{\eta}$,  is pre-quantizable.
Now a prequantum bundle $L$ is the bundle on $G/G_{\eta}$ determined by $\Lambda$.
So $\Gamma(L)$ can be identified with the space of $\Lambda$-equivariant functions on $G$.
In ${\cal O}_{\eta}$  one can consider loops $\psi$ of Hamiltonian
symplectomorphisms generated by  vector fields associated to elements of ${\frak g}$;
it turns out that  the  values of $\kappa$ on these
particular loops   
are related with the   character $\Lambda$. More precisely, 
if $\psi$ is the loop generated by  a family $A_t$ of elements of ${\frak g}$, 
from the action of $G$ on the $\Lambda$-equivariant functions  it is easy to 
prove that $\kappa(\psi)=\Lambda(h_1)$, where
$h_t$ is the
solution to Lax equation $\Dot h_th_t^{-1}=A_t$ (Theorem \ref{Tmcoador}). 
On the other hand, 
according to the Kostant version of Borel-Weil theorem \cite{bK66},
to each integral 
orbit ${\cal O}_{\eta}$ corresponds an irreducible representation $\pi$ of $G$,
whose highest weight is determined by $\eta.$ 
If $G$ is a semisimple Lie group, 
the choice of a maximal torus  contained in $G_{\eta}$
permits us to define a $G$-invariant complex structure $I$ on $G/G_{\eta}={\cal O}_{\eta}$.
This complex structure, in turn, determines a holomorphic structure on $L$, and  
  ${\cal Q}_I$ is just the space  $H^0(L)$ of holomorphic sections of $L$.  
When   $G_{\eta}$ is itself  a maximal torus,  
  the Borel-Weil theorem allows us to characterize the restriction of $\upsilon$ to
$H^0(L)$ in terms of its 
highest weight. It turns out that this restriction of $\upsilon$ is the representation $\pi^*$,
dual of  $\pi$.   It follows from Theorem \ref{Thchi1} that  the invariant $\kappa(\psi)$ for the closed isotopy
considered above is equal to $\chi(\pi^*)(h_1)/\text{dim}\,\pi$,
 where $\chi(\pi^*)$ is the character
of $\pi^*$ (Theorem \ref{BWWey}).    This result permits us 
to calculate $\kappa(\psi)$ using the Weyl's character formula. 

The paper is organized as follows.  
The second Section  is concerned with the new proof of the independence of 
$\kappa_x(\psi)$ from the point $x$. The existence of a lifting for   $\kappa$
to an ${\Bbb R}$-valued map is proved as well.
In Section 3 we study the map $\kappa$ when the manifold is a Hamiltonian $G$-space.
In the case that $G$ is a torus,   the existence of a moment map  such that its vertices
are integer lattice points is proved.
In Section 4 we study the invariant $\kappa(\psi)$  
 in an integral  coadjoint orbit of a compact Lie group $G$.       
Finally we apply some results of  this Section to two particular cases.
In the first one we determine the value of $\kappa(\psi)$ for a
closed isotopy $\psi$
in a general flag manifold, when the isotopy is
 generated by the action of the corresponding unitary group. 
  In \cite{aV01} we determined the value $\kappa(\psi)$ for
 a closed Hamiltonian flow $\psi$ in $S^2$, by direct calculation.
 Here we recover this number twice; by applying
Theorem \ref{Tmcoador}, and using Weyl's character  formula.    

 \smallskip

{\it Acknowledgements.} I thank the referee for his comments and for having pointed  me out
the references \cite{pS97} \cite{mS00}.


  \smallskip

\section{Action integral around  isotopies}

  We denote by $(M,\omega)$ a compact, connected,  symplectic, $2n$-dimensional manifold.
We will suppose that $(M,\omega)$ is {\it quantizable}, so   
 there exists
a Hermitian line bundle $\pi:L\rightarrow M$ (a prequantum bundle) with a connection,   
such that its curvature  
is $-2\pi i\omega$, that is, the first Chern class $c_1(L)=[\omega]$.

$L^{\times}=L-\{\text{zero section} \}$ is the corresponding ${\Bbb C}^{\times}$-principal
 bundle.
If $c\in{\Bbb C}=\text{Lie}({\Bbb C}^{\times})$ we denote by $W_c$ 
the vertical vector
field on $L^{\times}$ whose value at p is defined by the curve $\{p\cdot e^{2\pi itc}\}_t$.
If $f$ is a function on $M$ by $W_f$ is denoted the vector field on $L^{\times}$ given by
$W_f(p)=W_{f(\pi(p))}(p)$. As it is well-known
each section $\lambda$ of $L$ determines an equivariant  map 
$\lambda^{\sharp}:L^{\times}\rightarrow {\Bbb C}$ by the formula
$\lambda(\pi(p))=p\cdot\lambda^{\sharp}(p)$, and one has    
the known relations (see \cite{jS80} \cite{KN})
\begin{equation}\label{dericov}
 (D_X\lambda)^{\sharp}=X^{\sharp}(\lambda^{\sharp})=D_{X^{\sharp}}\lambda^{\sharp},\;\;\;
D_X\lambda=(\lambda^*\alpha)(X)\lambda,
\end{equation}
where $\alpha$ is the connection in the principal bundle $L^{\times}$.
It is also easy to prove that for $Y\in T_xM$
\begin{equation}\label{dericov1}
\lambda_*(Y)=Y^{\sharp}(\lambda(x))+W_c(\lambda(x)),
\end{equation}
with $c=(2\pi i)^{-1}(\lambda^*\alpha)(Y).$

On the other hand, if $X$ is the Hamiltonian vector field on $M$ determined by de function $f$
(i.e. $\iota_X\omega=-d\,f$), then the Lie derivative
${\cal L}_Z\alpha=0$, for $Z=X^{\sharp}-W_f$  (see \cite[page 56]{jS80}). 
 Hence $X^{\sharp}-W_f$ defines a $1$-parameter subgroup of diffeomorphisms of $L^{\times}$ which
preserve the connection $\alpha$.

\smallskip

Let $\{\psi_t\,|\, t\in[0,1] \}$ be the Hamiltonian isotopy associated to the time-dependent
Hamiltonian $f_t$.  This isotopy is generated by a
 unique  family of vector fields $X_t$ by the relations
$$\frac{d\psi_t(x)}{dt}=X_t(\psi_t(x)),\;\;\; \psi_0=\text{id},$$
with $\iota_{X_t}\omega=-d\,f_t.$
The Hamiltonian $f_t$ is determined by $\{\psi_t\}$ up to an additive constant. 
In the general case this constant will be fixed  
imposing that $\int_M f_t\omega^n=0$. 
The family $\{ X_t^{\sharp}-W_{f_t} \}_t$ of vector fields on $L^{\times}$
determines a set $\{F_t\}_t$ of preserving connection diffeomorphisms of $L^{\times}$,
such that $\pi\circ F_t=\psi_t\circ\pi.$ We call the isotopy $F_t$ the lift of   $\psi_t$.

Given $\tau$ a local section of $L^{\times}$, one defines $\tau_t:=F_t\circ\tau\circ\psi_t^{-1}.$
It is easy to prove the following relation between the corresponding equivariant maps
 \begin{equation} \label{equisharp}
\tau_t^{\sharp}=\tau^{\sharp}\circ F_t^{-1}.
\end{equation}

For each $t$ we have the map $\tau\in \Gamma(L^{\times})\mapsto \tau_t\in\Gamma(L^{\times})$.
This family of  maps can be regarded as a transport ${\cal T}_{\psi}$ along the isotopy $\{\psi_t\}$.

First of all we will determine the differential equation that defines 
 the transport ${\cal T}_{\psi}$.
Given the point $x\in M$ and the section $\tau$ of $L^{\times}$, 
 the tangent vector to the curve 
$\{F_t(\tau(x)) \}_t$ in $L^{\times}$ at the point $q_t=F_t(\tau(x))$ is, by definition of $F_t$,
\begin{equation}\label{tanvec}
Z_t(q_t):=X_t^{\sharp}(q_t)-W_{f_t}(q_t).
\end{equation}

On the other hand $F_t(\tau(x))=\tau_t(x_t)$, with  $x_t:=\psi_t(x)$;  hence 
 \begin{equation}\label{tanvec1}
Z_t(q_t)=\tau_{t*}(X_t(x_t))+W_a(q_t),
\end{equation}
with 
\begin{equation}\label{tauepsilon}
2\pi ia\cdot q_t= \frac{d\,\tau_{\epsilon}(x_t)}{d\,\epsilon} \biggr|_ {\epsilon=t} .
\end{equation}
The formula (\ref{dericov1}) applied to $\tau_t$ gives
\begin{equation}\label{tautXt}
\tau_{t*}(X_t(x))=X_t^{\sharp}(\tau_t(x))+W_c(\tau_t(x)),
\end{equation}
where $c=(2\pi i)^{-1}(\tau_t^*\alpha)(X_t(x_t))$.
It follows from (\ref{tanvec}), (\ref{tanvec1}) and (\ref{tautXt})
\begin{equation}\label{Wftqt}
-W_{f_t}(q_t)=W_c(q_t)+W_a(q_t).
\end{equation}
Using the identification of $T_{q_t}(L_{x_t})$ y $L_{x_t}$, where $L_{x_t}$ is
 the fiber of $L$ over $x_t$, one obtains from (\ref{Wftqt})
$$-2\pi if_t(x_t)\cdot q_t=  (\tau_t^*\alpha)(X_t(x_t))\cdot q_t+2\pi ia\cdot q_t.$$
By (\ref{dericov}) the first term on the right hand side is 
$(D_{X_t}\tau_t)(x_t)$. Hence by (\ref{tauepsilon})
$$\Big(\frac{d\,\tau_{\epsilon}}{d\,\epsilon}\Big)_{|\epsilon=t}(x_t)=-(D_{X_t}\tau_t)(x_t)-2\pi if_t(x_t)\cdot \tau_t(x_t).$$
As the point $x$ is arbitrary we have proved the following Proposition
that gives the differential equation for the transport ${\cal T}_{\psi}$ along $\psi$
\begin{Prop}\label{Profracd} The
    family of sections $\tau_t$ defined by $\tau_t=F_t\circ\tau\circ\psi_t^{-1}$ is 
the solution to 
  the equation 
\begin{equation}\label{paratrans}
\frac{d\,\tau_t}{d\,t}=-D_{X_t}\tau_t-2\pi if_t\tau_t,\;\;\; \tau_0=\tau.
\end{equation}
\end{Prop}

Let $X$ be the Hamiltonian vector field on $M$ associated to the function $f$. This vector field determines
the operator ${\cal P}_X:=-D_X-2\pi if$, which acts on the space $\Gamma(L)$.
So
$$\frac{d\,\tau_t}{dt}={\cal P}_{X_t}(\tau_t),\;\;\;\;\; \tau_0=\tau,$$
is the differential equation for the transport ${\cal T}_{\psi}$.
 
To prove the next Proposition we use the following Lemma; its proof is straightforward 

\begin{Lem}\label{L:LetH} Let $H$ be a diffeomorphism of $L$ over the diffeomorphism $\varphi$ of $M$.
If $H$ preserves the connection, then
$$ H_*(Y^{\sharp}(p))=\big( \varphi_*(Y) \big)^{\sharp}(H(p)),$$
for $Y\in T_xM$ and $p\in \pi^{-1}(x)$. And
$$H_*(W_f(p))=W_{f\circ\varphi^{-1}}(H(p)),$$
for any function $f$ defined on $M$.
\end{Lem}

\begin{Prop}\label{Avi8}
Let $\tau_t$ be the solution to (\ref{paratrans}), and $Y$ a vector field on $M$. Then
 $D_{Y_t}\tau_t=D_Y\tau$, for $Y_t=\psi_{t*}(Y)$.
\end{Prop}
\begin{pf}  
As $F_t$ preserves the connection, Lemma \ref{L:LetH} is applicable
$$Y_t^{\sharp}(\tau_t^{\sharp})=\big( \psi_{t*}(Y) \big)^{\sharp}(\tau_t^{\sharp})=
F_{t*}(Y^{\sharp})(\tau_t^{\sharp}).$$
Taking into account (\ref{equisharp})
$$Y_t^{\sharp}(\tau_t^{\sharp})=F_{t*}(Y^{\sharp})(\tau^{\sharp}\circ F_t^{-1})=
Y^{\sharp}(\tau^{\sharp}).$$
From (\ref{dericov}) it follows $D_{Y_t}\tau_t=D_Y\tau$. 
\end{pf}

Let $q$ be a point of $L^{\times}$ with $\pi(q)=x$, and $\tau$  a local section of $L^{\times}$
such that $x$ belongs to the domain of $\tau$. As $F_t(q)$ is in the 
fiber of $L$ over $x_t:=\psi_t(x)$,
then 
\begin{equation}\label{Ftq}
F_t(q)=\tau(x_t)\cdot m_t,
\end{equation}
 with $m_t\in{\Bbb C}^{\times}$. The curve $\{F_t(q)\}_t$ defines at the point
$q_t=F_t(q)$ the vector (\ref{tanvec}). On the other hand
 the curve $\tau(x_t)\cdot m_t$ defines also the vector
\begin{equation}\label{Rmt}
Z_t(q_t)=(R_{m_t})_*(\tau_*(X_t(x_t)))+W_c(q_t),
\end{equation}
where $c=(2\pi i)^{-1}(\Dot m_t/m_t)$ and $R_d$ 
is the right multiplication by $d\in{\Bbb C}^{\times}$ in $L^{\times}$. As $R_{d*}$ preserves 
the horizontal and vertical components,
by (\ref{dericov1})  
\begin{equation}\label{aux11}
Z_t(q_t)=X_t^{\sharp}(q_t)+W_d(q_t)+W_c(q_t),
\end{equation} 
with $d=(2\pi i)^{-1}(\tau^*\alpha)(X_t(x_t))$. It follows from 
the equality of  (\ref{aux11}) with (\ref{tanvec})
$$\frac{\Dot m_t}{m_t}=-(\tau^*\alpha)(X_t(x_t))-2\pi if_t(x_t).$$
That is,
\begin{equation}\label{mt}
m_t=\,\text{exp}\Big(-\int_0^t(\tau^*\alpha)(X_u(x_u))\,du-2\pi i\int_0^tf_u(x_u)\,du  \Big)
\end{equation} 

If the isotopy $\psi_t$ is closed, that is, $\psi_1=\text{id}$, then the curve $\psi_t(x)$
is nullhomologous  \cite[page 334]{Mc-S}.    Let $S$ be any 
 $2$-chain in $M$ whose boundary is the closed curve $\{\psi_t(x)\}_t$.
 As the curvature of the connection is $-2\pi i\omega$, by the Stokes theorem 
\begin{equation}\label{m1=m0}
m_1=\text{exp} \Big( 2\pi i\int_S\omega-2\pi i\int_0^1f_t(\psi_t(x))dt  \Big)
\end{equation}
 
 One obtains from   (\ref{Ftq}), (\ref{m1=m0})  and
the definition (\ref{kappa}) of $\kappa_x(\psi)$ 
\begin{equation}\label{F1q}
F_1(q)=q\cdot \kappa_{\pi(q)}(\psi).
\end{equation}
Therefore $F_1$ is the gauge transformation determined by the map
$\kappa_{-}(\psi):M\rightarrow U(1)$. As $F_1$ preserves the connection $\alpha$,
 this implies $d\,\kappa_{-}(\psi) =0$. Hence the map $\kappa_{-}(\psi)$ is
constant, and we have
\begin{Thm}\label{ThIn} If $\{\psi_t\}_{t\in[0,1]}$ is a closed Hamiltonian isotopy
in a quantizable manifold, the action
 integral along
the curve $\{\psi_t(x)\,|\, t\in[0,1]\}$ is independent of the point $x\in M$.  
\end{Thm}
As we said in Introduction a different proof
 of this result is given in \cite{pS97} (see also \cite[Section 3]{mS00}).

From (\ref{F1q}) and (\ref{equisharp}) it follows 
$\tau_1^{\sharp}=\kappa(\psi)\tau^{\sharp}$; that is,
\begin{Cor}\label{Corholo}
$\kappa(\psi)$ is the holonomy of the transport ${\cal T}_{\psi}$.
\end{Cor}

 As  we said ${\cal L}$ denote the space of loops in $\text{Ham}(M)$
based at $\text{id}$.
If $\psi,\,\xi\in{\cal L}$ and we denote by $\psi\ast \xi$ the usual product of paths, it
is immediate to check that 
\begin{equation}\label{pathproduct}
\kappa(\psi\ast\xi)= \kappa(\psi)\kappa(\xi). 
\end{equation}

Next we   will study the behavior of   $\kappa(\psi)$ 
under $C^{1}$-deformations of $\psi$.   We consider
the derivative of $\kappa(\psi^s)$ with respect to the parameter $s$ in a deformation
$\psi^s$ of $\psi$. That is, $\psi^s=\{\psi^s_t\,|\, t\in[0,1] \}$ is an isotopy with
$\psi^s_0=\psi^s_1=\text{id}$ generated by the time-dependent Hamiltonian $f^s_t$; 
furthermore we assume that
  $\psi^0=\psi$. By $\{X^s_t\}_t$ is denoted the family of Hamiltonian vector fields 
defined by $\{f^s_t\}_t$.

For $x\in M$ we put $\sigma^s(t):=\psi^s_t(x)$, so
 $\{\sigma^s(t)\,|\, t\in[0,1]\}$ is a closed curve and then
$$\kappa(\psi^s)=\text{exp}\,\Big(2\pi i\int_{S^s}\omega
  -2\pi i\int_0^1f_t^s(\sigma^s(t))\Big),$$              
where $S^s$ is a surface bounded by the curve $\sigma^s$. 
We set 
$$X_t:=X^0_t,\,\,\,  f_t:=f^0_t,\,\,\,  \sigma(t)=\sigma^0(t) .$$
The variation of $\sigma^s(t)$ with $s$
permits to define the vector fields $Y_t$; that is,
\begin{equation}\label{Yts}
Y_t(\sigma^s(t)):=\frac{\partial}{\partial s} \sigma^s(t).            
\end{equation}
For an ``infinitesimal" $s$ the curves $\sigma^l$, with $l\in[0,s]$ determine the 
``lateral surface" $J$ of one ``wedge"  whose base and cover
 are the surfaces $S$ and $S^s$ respectively. The ordered pairs of vectors   
$(X_t(\sigma(t)),Y_t(\sigma(t)))$
fix an orientation on $J$, which   in turn 
determines an orientation on the closed surface   $T=S\cup J\cup S^s$.
If we assume that $S$ and $S^s$ are oriented by means of the orientations 
of curves $\sigma$ and $\sigma^s$, from the fixed orientation on $T$ it 
 follows $T=J-S+S^s$. 

As $\omega$ satisfies the integrality condition
 \begin{equation}\label{derivadaS}
-\int_S\omega+\int_{S^s}\omega  =-\int_{J}\omega \,\,\text{(modulo}\,{\Bbb Z}\text{)}.
\end{equation}
Moreover 
\begin{equation}\label{intJOmega}
\int_{J}\omega=
s\int_0^1\omega\big(X_t(\sigma(t)),\,Y_t(\sigma(t)) \big)dt +O(s^2).
\end{equation}
On the other hand, for a given   $t\in[0,1]$  
\begin{equation}\label{derivfst}
\Big(\frac{d}{ds} f^s_t(\sigma^s(t))  \Big)_{|s=0}  =
\Big(\frac{\partial}{\partial s}f^s_t(\sigma(t))\Big)_{|s=0}+Y_t(\sigma(t))(f_t).
\end{equation}
We set
$$\Dot f_t(x):=\Big(\frac{\partial}{\partial s}f^s_t(x)\Big)_{|s=0}.$$
 As $\iota_{X_t}\omega=-df_t$, from
(\ref{derivfst}) it follows
\begin{equation}\label{d|dsint}
\frac{d}{ds}\Big|_{s=0} \int_0^1f^s_t(\sigma^s(t))dt = \int_0^1\Dot f_t(\sigma(t)) dt
-\int_0^1\omega\big(X_t(\sigma(t)),Y_t(\sigma(t))\big)dt.
\end{equation}
 By (\ref{derivadaS}), (\ref{intJOmega}) and  (\ref{d|dsint})
$$\kappa(\psi^s)-\kappa(\psi)=-2\pi is\kappa(\psi)\int_0^1\Dot f_t(\sigma(t)) dt+O(s^2),$$
 and finally
$$\Big(\frac{d}{ds}\kappa(\psi^s)\Big)_{|s=0}=
-2\pi i \kappa(\psi)\int_0^1\Dot f_t(\psi_t(x)) dt.$$
We have proved the following Proposition
\begin{Prop}\label{deformpsi}
 If $\psi^s$ is the deformation of $\psi$ defined by the family $f_t^s$ of 
time-dependent Hamiltonians, then
$$\frac{1}{2\pi i \kappa(\psi)}\Big(\frac{d}{ds}\kappa(\psi^s)\Big)_{|s=0}=
-\int_0^1\Dot f_t(\psi_t(x)) dt,$$
$x$ being any point of $M$.
\end{Prop}

\begin{Prop}\label{reprpi1}
$\kappa$ defines a group homomorphism
$\kappa:\pi_1(\text{Ham}(M))\rightarrow U(1)$.
\end{Prop}
\begin{pf}
With the above notations $\int_M f^s_t\omega^n=0$ for any $s$. Then 
$$\int_M(\Dot f_t\circ\psi_t)\omega^n=0.$$
 By Proposition \ref{deformpsi}
$$\int_0^1\Dot f_t(\psi_t(x)) dt$$
 is independent of the point $x\in M$. Hence
$$\Big(\frac{d}{ds}\kappa(\psi^s)\Big)_{|s=0}\int_M\omega^n=
-2\pi i\kappa(\psi)\int_0^1 dt\int_M(\Dot f_t\circ\psi_t)\omega^n=0.$$
So $\kappa(\psi)$ depends only on the homotopy class $[\psi]\in\pi_1(\text{Ham}(M))$.
From (\ref{pathproduct}) we conclude that $\kappa$ is a group homomorphism.
\end{pf}

As a consequence of Proposition \ref{reprpi1} it makes sense to define the {\it action integral}
around an element $[\psi]\in\pi_1(\text{Ham}(M))$ as $\kappa_x(\psi)$, for $x$ an arbitrary point of $M$.

\smallskip


\section {Hamiltonian $G$-spaces.}

Let $G$ be a compact, connected Lie group which acts on the left on the quantizable 
manifold  
  $(M,\omega)$ by symplectomorphisms.   We assume that this action is
Hamiltonian, and that $\Phi:M\rightarrow {\frak g}^*$ is a   moment  
 map  for this action. That is, $M$ is a Hamiltonian $G$-space \cite{vG94}.

Given $A\in{\frak g}$, we denote by $X_A$,  the vector field 
on $M$ generated by $A$. Then  $\langle d\Phi(Y),\, A\rangle=\omega(Y, X_A)$, 
for any vector field $Y$ on $M$. The $A$-component of $\Phi$    will be denoted by $f_A$;
that is, $f_A(x)=\langle \Phi(x),\,A \rangle$. Hence
\begin{equation}\label{moment}
\iota_{X_A}\omega=-d\,f_A \,\,\,\,\text{and}\,\,\,\, \{f_A,f_B\}:=\omega(X_ B,X_ A)=f_{[A,B]}.
\end{equation}

Now the prequantization map ${\cal P}$ restricted to ${\frak g}$ is 
\begin{equation}\label{preqmap}
 A\in{\frak g}\mapsto {\cal P}_A=-D_{X_ A}-2\pi if_A\in \text{End}(\Gamma(L)),
\end{equation}
  and using  (\ref{moment}) it is straightforward to prove the following Proposition
\begin{Prop}
The map ${\cal P}$ is a Lie algebra homomorphism.
\end{Prop}

Given a family $\{g_t\}_{t\in[0,1]}$ of elements of $G$ with $g_0=e$, it
 determines a Hamiltonian isotopy $\varphi_t:M\rightarrow M$
by 
\begin{equation}\label{g_tx}
\varphi_t(x)=g_t\cdot x.
\end{equation}
It is easy to prove the following Proposition that gives the vector fields
determined by the isotopy $\varphi_t$.
\begin{Prop}\label{P1}
 The isotopy  (\ref{g_tx}) is defined by the equations
$$ \frac{d  \varphi_t}{dt}=X_{A_t}\circ\varphi_t,\;\;\; \varphi_0=\text{id}, $$
where $A_t$ is the element of ${\frak g}$ given by $A_t:=\Dot g_tg_t^{-1}$. 
\end{Prop}

If $x_0\in M$ is a fixed point for the $G$-action, then $\Phi(x_0)$ is called a 
{\em vertex} of $\Phi$
\cite{vG820}. If $\varphi$ is a closed isotopy we will express  
 the $U(1)$-valued action integral $\kappa(\varphi)$   in terms of the vertices of $\Phi$.

\begin{Thm} \label{Th1} Let $(M,\omega)$ be a compact, quantizable $G$-space. If $x_0$ is a fixed point for the $G$-action and if   
 the isotopy $\varphi_t$ defined by
(\ref{g_tx}) is closed (i.e. $\varphi_1=\text{id}$), then
$$\kappa(\varphi)=\text{exp}\big(-2\pi i\langle \Phi(x_0),\,\int_0^1A_t dt\rangle\big),$$
with $A_t=\Dot g_t g_t^{-1}$
\end{Thm}
\begin{pf} In this case the curve $\{\varphi_t(x_0)\}$ is a point, so integral of $\omega$ in
(\ref{kappa}) vanishes. On the other hand,  by Proposition \ref{P1}
 the corresponding   time-dependent Hamiltonian is $f_{A_t}=\langle \Phi,\,A_t\rangle$.

\end{pf}

Now we assume that $G=T$ is the $r$-torus.
\begin{Thm} \label{Themfet} Let $(M,\omega)$ be a compact, quantizable manifold
equipped with a Hamiltonian action of the $r$-torus $T$.
If $v$ and $v'$ are vertices of a moment map $\Phi$ for 
the $T$-action, then $v-v'$ belongs to the lattice ${\Bbb Z}^r$ of ${\frak t}^*$.
\end{Thm}
\begin{pf} Let $x_0$ and $x_1$ be fixed points with $\Phi(x_0)=v$ and $\Phi(x_1)=v'$. If $A$ is a 
 vector  
of ${\frak t}^*$ with integer coordinates, then 
$$\{\phi_t(x)=e^{2\pi it A}\cdot x  \}_{t\in[0,1]}$$
is a closed isotopy of $(M,\omega)$. Hence, by the independence of $\kappa_x(\phi)$
 from the point $x$ 
and Theorem \ref{Th1}, we conclude
$$\text{exp}\big(4\pi^2 \langle \Phi(x_0),\,A \rangle    \big)=
\text{exp}\big(4\pi^2\langle \Phi(x_1),\,A \rangle    \big).$$
That is, $\langle v-v',\,A \rangle =0$. This relation holds
for any $A\in{\frak t}^*\cap{\Bbb Z}^r$.
Hence  $v-v'\in{\Bbb Z}^r$. 

 \end{pf}

 \begin{Cor}\label{Corigene} Under the hypotheses of Theorem \ref{Themfet},
  there is a moment map $\Phi$ for
the $T$-action such that the vertices of $\Phi$ are integer lattice points. 
\end{Cor} 
\begin{pf} If $\tilde\Phi$ is a moment map and $x_0$ is a fixed point, we set
 $v:=\tilde\Phi(x_0)$. We define the map $\Phi:=\tilde\Phi-v$, which is also a moment map,
since $T$ is Abelian. Hence $0$ is a
vertex of $\Phi$. By Theorem \ref{Themfet}, for any vertex $v'$ of $\Phi$,
$v'=v'-0\in{\Bbb Z}^r$
\end{pf}

In particular, if the toric manifold associated to a Delzant polytope $\Delta$
is quantizable, then the vertices of $\Delta$ are integer lattice points. So we recover a 
well-known  result relative to Delzant spaces 
  \cite{tD88} \cite{vG94}.

\smallskip

The action of $G$ on $M$    is  said to be  {\em pre-quantizable} 
  if there is a global action of $G$
on the prequantum bundle $L$ which induces the action (\ref{preqmap}) of ${\frak g}$ 
on the space $\Gamma(L)$   \cite{vG820} (In  \cite{vG82} one says that the prequantum data
$(L,D)$ are $G$-invariant).
The thesis of Corollary \ref{Corigene} has been proved in
\cite[$\S$ 8, Corollary 1]{vG820}
  under the additional assumption that the
$G$-action is pre-quantizable.

  Henceforth in this Section we assume the existence of this lift of the $G$-action to $L$. We denote
by $\rho(g)$ (resp. $\varphi_g$) the diffeomorphism of $L$ (resp. $M$) associated to $g\in G$.
So 
$\pi\circ\rho(g)=\varphi_g\circ\pi$.
We denote by $\upsilon$ the representation of $G$ on $\Gamma(L)$ induced by $\rho$. That is,
for $g\in G$ and $\tau\in\Gamma(L)$
\begin{equation}\label{upsdef}
\upsilon(g)(\tau)=\rho(g)\circ\tau\circ\varphi_g^{-1}.
\end{equation}
 Moreover, given a curve  $\{g_t\}_{t \in[0,1]}$ in $G$ starting at $e$,  
    the fact that $\rho$ induces ${\cal P}$ implies that
\begin{equation}\label{rho(gt)}
\frac{d}{dt}\biggr|_ {t=0}\upsilon(g_t)\tau   ={\cal P}_A\tau,
\end{equation}
where $A\in{\frak g}$  is the derivative $\Dot g_t$ at $t=0$.

As we said in Section 2, the Hamiltonian isotopy $\varphi_t=\varphi_{g_t}$
admits a lift to     an isotopy 
$F_t$ of $L^{\times}$. By Proposition \ref{P1} the family of vector fields
$$Z_t:=X_{A_t}^{\sharp}-W_{f_t},$$  
where $A_t:=\Dot g_tg_t^{-1}$ and $f_t=\langle\Phi,\,A_t\rangle$, generates the 
isotopy $F_t$.
 
It is proper to ask if the diffeomorphisms $F_t$ and $\rho(g_t)$ are related. We will
prove that $F_t=\rho(g_t)$.

\begin{Lem}\label{Lemasal}
Let  $h_t$ be a curve in $G$ with $h_0=e$ and $\Dot h(0)=B\in{\frak g}$, and let $p$ be
a point in $L^{\times}$. Then 
$$\frac{d\,\rho(h_t)(p)}{dt}\biggr|_{t=0}= X^{\sharp}_B(p)-W_{f_B}(p).$$  
\end{Lem}
\begin{pf}
We will determine the vector tangent to $\{\rho(h_t)p\}_t$ at $t=0$.
Let $\tau$ be a section such that $\tau(x)=p$, then 
\begin{equation}\label{kaaux}
\frac{d}{dt}\,\rho(h_t)\tau(h_t^{-1}x)\biggr|_{t=0}=
\frac{d}{dt}\,\rho(h_t)p  \biggr|_{t=0}+\tau_*(-X_B(x)).
\end{equation}
It follows from (\ref{dericov1})  
\begin{equation}\label{RemarkSect2}
\tau_*(X_B(x))=X_B^{\sharp}(\tau(x))+(\tau^*\alpha)(X_B)\tau(x)=
X_B^{\sharp}(\tau(x))+(D_{X_B}\tau)(x).
\end{equation}
The left hand side in (\ref{kaaux}) is equal to $({\cal P}_B\tau)(x)$ by (\ref{rho(gt)})
and (\ref{upsdef}).
Therefore (\ref{kaaux}) and (\ref{RemarkSect2}) give rise to 
$$\frac{d\,\rho(h_t)(p)}{dt}\biggr|_{t=0}=X_{B}^{\sharp}(p)-W_{f_B}(p).$$
 \end{pf}

\begin{Prop}\label{P:If} The isotopy $\rho(g_t)$ is generated by the family 
of vector fields $Z_t$. So $\rho(g_t)=F_t$.
\end{Prop}
\begin{pf}
Given $t$, we put $h_{t'}:=g_{t'}g_{t}^{-1}$; then $\Dot h_t=A_t$. Since $\rho$ is 
a representation of $G$, by Lemma \ref{Lemasal} 
$$\frac{d\,\rho(g_{t'})}{dt'}\biggr|_{t'=t}(\rho(g^{-1}_t(p))=
\frac{d\,\rho(h_{t'})(p)}{dt'}\biggr|_{t'=t}= X_{A_t}^{\sharp}(p)-W_{f_t}(p).$$
That is, 
$$\frac{d\,\rho(g_t)}{dt}=Z_t\circ \rho(g_t).$$
\end{pf}

 The following Corollary asserts that the transport along $\varphi$
is defined by the representation $\upsilon$.

\begin{Cor}\label{Cori1} If $g_t$ is a family of elements in
 $G$ with $g_0=e$ and $A_t:=\dot g_tg_t^{-1}$,
then the solution to
\begin{equation}\label{ransG}
\frac{d\,\tau_t}{dt}={\cal P}_{A_t}(\tau_t),\;\;\; \tau_0=\tau
\end{equation}
  is $\tau_t=\upsilon(g_t)(\tau).$  In particular, $\tau_1$ depends only on the endpoint $g_1$
of the curve $g_t$.
\end{Cor}
\begin{pf} 
 The statement follows from 
Proposition \ref{Profracd} together with Proposition \ref{P:If} and (\ref{upsdef}).
\end{pf}

\smallskip

Now we will express the $U(1)$-action integral $\kappa(\varphi)$ in terms 
of the representation $\upsilon$.
We will construct finite-dimensional $\upsilon$-invariant subspaces of $\Gamma(L)$, and
$\kappa(\varphi)$ will be related with the characters of the restrictions of $\upsilon$ to
these subspaces.

An almost complex structure $J$ on $M$ is called {\it compatible} with
$\omega$ if $\omega(J.\,,\,J.\,)=\omega(\,.\, ,.\,)$ and $\omega(\,.\,,J.\,)$ is
positive definite. Since  
 $G$ is a compact group one can construct $G$-invariant compatible almost complex structures
on $M$. Let I be such an   almost complex structure. A section $\tau$   of $L$
is said to be $I$-polarized if
$D_X\tau=0$, for any vector field $X$ on $M$ of type $(0,1)$ relative to $I$. We set   
${\cal Q}_I$  for the space of $I$-polarized sections.  

With $\Omega^{0,k}(L)$ we denote the space of $(0,k)$-forms on $M$ with  values in $L$. 
The operator 
$\bar D=(1+iI)D$ extends in the usual manner to an operator
$$\bar D:\Omega^{0,k}(L)\rightarrow \Omega^{0,k+1}(L).$$
Although $(\Omega^{0,*}(L),\,\bar D)$ is not a complex,
the Riemannian metric $\omega(\,.\,,\,I.\,)$ on $M$ and the Hermitian metric 
on $L$ determine the adjoint operator $\bar D^*$ of $\bar D$; moreover the  
operator
$$\bar D+\bar{D}^*:\Omega^{0,\text{even}}(L)\rightarrow \Omega^{0,\text{odd}}(L)$$
is elliptic (see \cite[page 75]{vG94}). Since  
$${\cal Q}_I=\text{Ker}\big((\bar{D}+\bar D^*)_{|\Omega^{0,0}}\big),$$
 then ${\cal Q}_I$ is a finite dimensional vector space.

Given $g\in G$, $g$ is the endpoint of a curve $g_t$ in $G$ with $g_0=e$.
By Corollary \ref{Cori1} $\upsilon(g)(\tau)=\tau_1$, where $\tau_t$ is the solution of 
(\ref{ransG}).
If $Y$ is a vector  of type $(0,1)$ with respect to $I$, by the $G$-invariance of $I$
the vector 
$\varphi_{t*}(Y)$ is  of type $(0,1)$ as well,   $\varphi_t$ being the diffeomorphism 
defined  in (\ref{g_tx}).
 This fact together with Proposition \ref{Avi8} imply that 
$\tau_1\in{\cal Q}_I$, if $\tau$ is $I$-polarized. Hence ${\cal Q}_I$
 is a $\upsilon$-invariant subspace of $\Gamma(L)$.
We set $\upsilon_I$ for the restriction of $\upsilon$ to ${\cal Q}_I$.

If the isotopy defined  in   (\ref{g_tx}) is closed,  
 from  Corollary  \ref{Corholo} and Corollary \ref{Cori1} it
follows   
$\upsilon_I(g_1)(\tau) =\kappa(\varphi)\tau.$
   Hence the character $\chi(\upsilon_I)$ 
at the point $g_1$ equals  $\kappa(\varphi)\,\text{dim}\,{\cal Q}_I$. 
One has the following Theorem

\begin{Thm}\label{Thchi1}  If $M$ is a Hamiltonian $G$-space and the $G$-action is
pre-quantizable, then  the action integral $\kappa(\varphi)$
 around the closed isotopy
$\varphi_t(x)=g_t\cdot x$ is equal to 
$$\frac{\chi(\upsilon_I)(g_1)}{ \text{dim}\,{\cal Q}_I },$$ 
where $I$ is any $G$-invariant, compatible almost complex structure on $(M,\omega)$.
\end{Thm}



\bigskip

\section {The invariant $\kappa(\psi)$ in a coadjoint orbit} 

\medskip
 
Let $G$ be a compact Lie group, and we consider the coadjoint action of $G$
 on ${\frak g}^*$ defined by 
$$(g\cdot\eta)(A)=\eta(g^{-1}\cdot A),$$
for $g\in G$, $\eta\in{\frak g}^*$, $A\in{\frak g}$ and $g\cdot A=\text{Ad}_g A$
(see \cite{aK76} \cite{nW92}).

If $X_A$ is the vector field on ${\frak g}^*$ determined by $A$, the map
 $l_g:\mu\in{\frak g}^*\mapsto g\cdot \mu\in{\frak g}^*$
satisfies
\begin{equation}\label{l_g_*}
(l_g)_*(X_A(\mu))=X_{g\cdot A}(g\cdot \mu).
\end{equation}

Given $\eta\in{\frak g}^*$, by ${\cal O}_{\eta}=:{\cal O}$ will be denoted the orbit of $\eta$
under the coadjoint action of $G$.   
 On ${\cal O}$
one can consider the $2$-form $\omega$ determined by
\begin{equation}\label{omeganu}
 \omega_{\nu}(X_{A}(\nu), X_{B}(\nu))=\nu([A,B]).
\end{equation}
This $2$-form determines a symplectic structure on ${\cal O}$, and the action of $G$
preserves $\omega$. For each $A\in{\frak g}$ one defines 
the function $h_A\in C^{\infty}({\cal O})$ by $h_A(\nu)=\nu(A)$, and for this function holds the 
formula
\begin{equation}\label{iotah}
\iota_{X_{A}}\omega=dh_A.
\end{equation}

The orbit ${\cal O}$ can be identified with $G/G_{\eta}$, where $G_{\eta}$
is the subgroup of isotropy of
 $\eta$. The Lie algebra of this subgroup is    
$${\frak g}_{\eta}=
\{A\in{\frak g}\,|\, \eta([A,B])=0, \,\,\text{for every}\,\,B\in{\frak g}  \}$$

One says that the linear functional
\begin{equation}\label{intfunc}
\lambda:C\in{\frak g}_{\eta}\mapsto 2\pi i\eta(C)\in i{\Bbb R}
\end{equation}
is integral iff there is a character $\Lambda:G_{\eta}\rightarrow U(1)$
whose derivative is the functional (\ref{intfunc}) (see \cite{bK70}).
Henceforth we assume the existence of such a character $\Lambda$.
We will prove that  
the orbit  ${\cal O}$ possesses a $G$-invariant prequantization. 
A prequantum bundle $L$ over ${\cal O}=G/G_{\eta}$ is defined by
$L=G\times_{\Lambda}{\Bbb C}=(G\times {\Bbb C})/\simeq$, with
$(g,z)\simeq (gb^{-1}, \Lambda(b)z)$, for $b\in G_{\eta}$.

Each section $\sigma$ of $L$ determines a $\Lambda$-equivariant function $s:G\rightarrow{\Bbb C}$
by the relation
\begin{equation}\label{sfunction}
\sigma(gG_{\eta})=[g,s(g)].
\end{equation}

The ${\Bbb C}^{\times}$-principal bundle associated to 
$L$ is $L^{\times}=L-\{\text{zero section}\}$. 
The lift $\sigma^{\sharp}:L^{\times}\rightarrow {\Bbb C}$ 
of the section $\sigma$ and its corresponding 
  $\Lambda$-equivariant function $s$  
are related by the formula
\begin{equation}\label{ssostenido}
s(g)=\sigma^{\sharp}([g,z])z.
\end{equation}

If $v$ denotes the element $[e,1]\in L^{\times}$, then 
$T_v(L^{\times})\simeq({\frak g}\oplus {\Bbb C})/{\frak f}_v$,
with
$${\frak f}_v=\{(B,\,-2\pi i\eta(B))\,|\, B\in{\frak g}_{\eta}  \}.$$
The connection form $\Omega$ on $L^{\times}$ is constructed in \cite{bK70} p.198.  
The form $\Omega$ can be written $\Omega=(\theta, \,d)$, where $\theta$ is the left 
invariant form on $G$ whose value at $e$ is $\eta$, and 
$d\in\text{Hom}_{\Bbb C}({\Bbb C},{\Bbb C})$ is defined by $d(z)=(2\pi i)^{-1}z$.
 It is clear that 
$\Omega_v$ vanishes on ${\frak f}_v$ and that it defines an element of $T_v^*(L^{\times})$.

 We denote by ${\cal E}_{\Lambda}$ the space of $\Lambda$-equivariant functions on $G$.
The identification $\Gamma(L)\thicksim {\cal E}_{\Lambda}$ allows us to translate the 
action ${\cal P}$ defined in (\ref{preqmap}) to a
 representation of ${\frak g}$ on ${\cal E}_{\Lambda}$.

\begin{Thm}\label{thms} The action ${\cal P}$ on ${\cal E}_{\Lambda}$ 
 is given by ${\cal P}_A(s)=-R_A(s)$,
where $R_A$ is the right invariant vector field on $G$ determined by $A$.
\end{Thm}
\begin{pf}
Let $\sigma$ be a section of $L$, by (\ref{iotah}) 
${\cal P}_A(\sigma)=-D_{X_A}\sigma+2\pi ih_A\sigma.$ 
We will determine the lift $({\cal P}_A(\sigma))^{\sharp}$.

 The vector $X_{A}(g\cdot\eta)\in T_{g\cdot\eta}({\cal O})$ is defined by the curve
$u\mapsto e^{uA}g\cdot\eta$ in ${\cal O}$. A lift of this curve at the point $[g,z]\in L^{\times}$
will be a curve of the form
$\gamma(u)=[e^{uA}g,\, z_u]$, with $z_u=ze^{ux}$. The vector tangent to $\gamma$ at $[g,z]$ 
is $\Dot\gamma(0)=[R_{A}(g),x]$, where $R_{A}(g)$ is 
the value at $g$ of the right invariant vector field in $G$ defined by $A$.

The condition $\Omega(\Dot\gamma(0))=0$ implies 
\begin{equation}\label{x0-2}
x=-2\pi i\eta(g^{-1}\cdot A).
\end{equation}
Therefore the horizontal lift of $X_{A}(g\cdot\eta)$ is
$$X_{A}^{\sharp}([g,z])=[R_{A}(g),\, -2\pi i\eta(g^{-1}\cdot A)],$$
and by  (\ref{ssostenido}) the action of $X_{A}^{\sharp}([g,z])$ on 
the function $\sigma^{\sharp}$
can expressed in terms of $s$
$$X_{A}^{\sharp}([g,z])(\sigma^{\sharp})=
\frac{d}{du}\biggr|_{u=0}\Big( \frac{s(e^{uA}g)}{ze^{ux}} \Big)=\frac{R_{A}(g)(s)}{z}-
\frac{xs(g)}{z}.$$
Since $X_{A}^{\sharp}(\sigma^{\sharp})=(D_{X_A}\sigma)^{\sharp}$, from
(\ref{x0-2})
and (\ref{ssostenido}) it turns out that 
 the equivariant function associated to $D_{X_{A}}\sigma$ is 
\begin{equation}\label{equivarDat}
g\in G\mapsto R_{A}(g)(s)+2\pi i\eta(g^{-1}\cdot A)s(g)\in{\Bbb C}.
\end{equation}
 
Obviously the equivariant function defined by the section $h_{A}\sigma$ is
the function $\lambda_{A}s$, where
 $\lambda_{A}(g)=h_{A}(gG_{\eta})=(g\cdot\eta)(A)=\eta(g^{-1}\cdot A).$
It follows from (\ref{equivarDat})   that
the equivariant function which corresponds to
$-D_{X_A}\sigma+2\pi ih_A\sigma$ is $-R_A(s)$.
\end{pf}

\begin{Cor}\label{CorII} The action ${\cal P}$ on ${\cal E}_{\Lambda}$ is induced by the action
$$\upsilon:(b,s)\in G\times {\cal E}_{\Lambda}\mapsto s\circ{\cal L}_{b^{-1}}\in{\cal E}_{\Lambda},$$
where ${\cal L}_c$ is left multiplication by $c$   in the group $G$.
 \end{Cor}
\begin{pf}If $g_t=e^{tA}\in G$, then
$$\frac{d\,\upsilon_{g_t}(s)}{dt}\biggr|_{t=0}(g)=
\frac{d}{dt}\biggr|_{t=0}s(e^{-tA}g)=-R_A(g)(s)={\cal P}_A(s)(g).$$
\end{pf} 

\begin{Cor}\label{CorIII} 
The action of $G$ on ${\cal O}_{\eta}$ is pre-quantizable.
\end{Cor}

\begin{pf}
On $L$ we define the following representation of $G$,  
$\rho(g')([g,z])=[g'g,z]$. For $g_t=e^{tA}$ and the section $\sigma$ of $L$
$$\frac{d}{dt}\biggr|_{t=0}\rho(g_t)\sigma(g_t^{-1}gG_{\eta})=
\frac{d\,[g,\,s(e^{-tA}g)]}{dt}\biggr|_{t=0}  =[g,\,-R_A(g)(s)]={\cal P}_A(\sigma)(gG_{\eta}).$$
The Corollary follows from (\ref{upsdef}) and (\ref{rho(gt)}).
\end{pf}

\smallskip

Let $\{\psi_t\,|\, t\in[0,1] \}$ be a closed Hamiltonian isotopy on ${\cal O}$.
 We also assume that
 the corresponding Hamiltonian vector fields are generated by elements of ${\frak g}$;
 that is,
$$\frac{d\psi_t(q)}{dt}=X_{A_t}(\psi_t(q)),\,\,\,\text{with}\,\,\,A_t\in{\frak g}.$$
If $\sigma$ is a section of $L$,  $\sigma_t$ will denote the solution to the equation
\begin{equation}\label{transportt}
\frac{d\sigma_t}{dt}= {\cal P}_{A_t}(\sigma_t),\,\,\,\,\,\,\, \sigma_0=\sigma.
\end{equation}
 By Theorem \ref{thms}, equation (\ref{transportt}) on the points
$\{h_t\}_{t\in[0,1]}$ of a curve in $G$ gives rise to   
\begin{equation}\label{Dotst}
\Dot s_t(h_t)=-R_{A_t}(h_t)(s_t),
\end{equation}
for the corresponding equivariant functions.
In particular, if $h_t$ is the curve such that $h_0=e$ and 
$\Dot h_t=R_{A_t}(h_t)\in T_{h_t}(G);$ in other words, $h_t$ satisfies the Lax equation
    $\Dot h_th_t^{-1}=A_t$, then
$$R_{A_t}(h_t)(s_t)=\frac{d}{du}\biggr|_{u=t}s_t(h_u).$$
Using (\ref{Dotst}) one deduces
\begin{equation}\label{constanc}
\Dot s_t(h_t)+\Dot h_t(s_t)=0
\end{equation}
If we consider the function $w:[0,1]\rightarrow{\Bbb C}$ defined by
$w_t=s_t(h_t)$; by
(\ref{constanc}) $w$ is constant. So $s_1(h_1)=s_0(e).$
If $h_1\in G_{\eta}$, as $s_1$ is $\Lambda$-equivariant
$s_1(h_1)=\Lambda(h_1^{-1})s_1(e)$; so
\begin{equation}\label{sigma1}
\sigma_1(eG_{\eta})=\Lambda(h_1)\sigma_0(eG_{\eta}).
\end{equation} 
The following Theorem, which gives the invariant $\kappa(\psi)$
in terms of $\Lambda$,  is consequence of Corollary \ref{Corholo} and (\ref{sigma1})
\begin{Thm}\label{Tmcoador} If $\{\psi_t \}$ is the closed 
Hamiltonian isotopy in ${\cal O}$ generated by the 
vector fields $\{X_ {A_t} \}$, then $\kappa(\psi)=\Lambda(h_1)$, where $h_t\in G$
is the solution to $\Dot h_th_t^{-1}=A_t$, with $h_0=e$ and $h_1\in G_{\eta}$.
\end{Thm} 

\smallskip

Let us assume that $G$ is semisimple Lie group \cite{wF91}, and let $T$ a maximal torus
with $T\subset G_{\eta}$ (see \cite{vG84} p.166). One has the standard 
decomposition of ${\frak g}_{\Bbb C}={\frak g}\otimes_{\Bbb R}{\Bbb C}$
in direct sum of root spaces
$${\frak g}_{\Bbb C}={\frak h}\oplus\sum {\frak g}_{\alpha},$$
where ${\frak h}={\frak t}_{\Bbb C}$, and $\alpha$ ranges over the set of roots.

 We denote by $\alpha^{\vee}$ the element
of $[ {\frak g}_{\alpha} , {\frak g}_{-\alpha}]$ such that $\alpha(\alpha^{\vee})=2$.
On the other hand $\eta$ extends in a natural way to ${\frak g}_{\Bbb C}$, and 
 if $Y\in{\frak g}_{\alpha}$, then 
$$0=\eta([\alpha^{\vee},Y])=2\eta(Y).$$
Hence $\eta$ vanishes on $\sum {\frak g}_{\alpha}.$ If $\eta(\alpha^{\vee})\ne 0$, for all
root $\alpha$, then ${\frak g}_{\eta}={\frak t}$; in this case $\eta$ is said to be regular.
Henceforth we assume that $\eta$ is regular. Let $P$ be the set of roots $\alpha$ such
that $\eta(\alpha^{\vee})<0$. Then the real counterpart of the
above direct sum decomposition is
$${\frak g}={\frak t}\oplus\sum_{\alpha\in P}
\big({\frak g}_{\alpha}\oplus {\frak g}_{-\alpha}\big)\cap{\frak g}.$$

We define ${\frak b}={\frak h}\oplus {\frak n}$, where 
$${\frak n}=\sum_{\alpha\in P}{\frak g}_{\alpha}.$$
Then ${\frak b}$ is a Borel subalgebra of ${\frak g}_{\Bbb C}$, which corresponds
to a Borel subgroup $B$ of $G$.

We have 
$$T_{\eta}({\cal O})={\frak g}/{\frak g}_{\eta}=
\sum_{\alpha\in P}\Big( {\frak g}_{\alpha}\oplus{\frak g}_{-\alpha} \Big)\cap{\frak g}.$$
Hence 
$$T^{\Bbb C}_{\eta}({\cal O})=\sum_{\alpha\in P}\Big(  {\frak g}_{\alpha}\oplus{\frak g}_{-\alpha} \Big ).$$
One defines 
$$T^{0,1}_{\eta}{\cal O}:={\frak n},$$
and 
$$T^{0,1}_{g\cdot\eta}{\cal O}:=\{ X_{g\cdot A}(g\cdot\eta)\,|\,A\in{\frak n}\}.$$
If $g_1\cdot\eta=g_2\cdot \eta$, then $g_1^{-1}g_2\in T$.
As ${\frak g}_{\alpha}$ is an eigenspace for the action of $T$, then 
$g_1^{-1}g_2\cdot A\in{\frak n}$, if $A\in{\frak n}$. Therefore the
spaces $T^{0,1}_{g\cdot\eta}$ are well-defined.

 For  
$A\in{\frak n}$, one can define the vector field ${\cal A}$ on  ${\cal O}$
by ${\cal A}(g\cdot\eta)=X_{g\cdot A}(g\cdot \eta)$.
By (\ref{l_g_*}) $(l_g)_*{\cal A}={\cal A}$, hence the above complex foliation
defined on ${\cal O}$ is $G$-invariant. Since
the vector $X_{g\cdot A}(g\cdot\eta)$ is defined by the curve
$e^{tg\cdot A}g\cdot \eta=ge^{tA}\cdot \eta$, then     the left invariant vector field $L_A$
on $G/T$ is the field which
 corresponds to  ${\cal A}$,
  in the identification of $G/T$ with ${\cal O}$.

The vector spaces $T^{1,0}$ are defined in the obvious way. As ${\frak n}$
is a subalgebra of ${\frak g}_{\Bbb C}$, the decomposition
$T^{\Bbb C}({\cal O})= T^{1,0}\oplus T^{0,1}$ define a complex structure $I$ on 
${\cal O}$. This complex manifold can be identified with $G_{\Bbb C}/B$.
  
Using the complex structure on ${\cal O}=G/T$ and the covariant derivative $D$
 on the prequantum bundle $L=G\times_{\Lambda}{\Bbb C}$, it is possible to
define a holomorphic structure in $L$. The section $\tau$ of $L$
is said to be holomorphic iff $D_Z\tau=0$ for any vector field $Z$
of type $(0,1)$. In this way $L$ can be regarded as a holomorphic line bundle
over $G_{\Bbb C}/B$, and 
with the notation of Section 3 ${\cal Q} _I=H^0(G_{\Bbb C}/B, L)$.

The homomorphism $\Lambda:T\rightarrow U(1)$
extends trivially to $B$, since $B$ is a semidirect product of $H=T_{\Bbb C}$ and
the nilpotent subgroup whose Lie algebra is ${\frak n}$. And  each   section $\sigma$
of $L$ determines a  function $s:G_{\Bbb C}\rightarrow {\Bbb C}$ which
is $\Lambda$-equivariant. On the other hand,
given $A\in {\frak n}$, the Proof of Theorem \ref{thms} shows that 
the equivariant function associated
to $D_{\cal A}\sigma$ is the map
$$g\in G_{\Bbb C}\mapsto R_{g\cdot A}(g)(s)+2\pi i\eta(g^{-1}g\cdot A)s(g)\in{\Bbb C}.$$
As $\eta$ vanishes on $\frak n$ and the vectors $R_{g\cdot A}(g)$ and $L_A(g)$
 are equal, the equivariant function associated
to $D_{\cal A}\sigma$ is $L_A(s)$. 
Therefore if $\sigma$ is holomorphic, then $L_A(s)=0$
for any $A\in{\frak n}$; that is, 
 $s$ is a holomorphic function on $G_{\Bbb C}$. 
 So 
  the space $H^0(G_{\Bbb C}/B, L)$ is isomorphic to the space
$${\cal E}_{\Lambda,I}:=\{s:G_{\Bbb C}\rightarrow {\Bbb C}\,|
\, s \,\,\text{is holomorphic and}\,\,\,\Lambda-\text{equivariant}   \}.$$

The Borel-Weil Theorem asserts that the action of $G$ on 
the space ${\cal E}_{\Lambda,I}$ given by 
$g\star s=s\circ{\cal L}_{g^{-1}}$ is an irreducible representation of $G$;
more precisely the contragredient representation of that one whose highest weight
is $-\lambda $ (see \cite{JJ00} pages 290, 300).

Denoting by $\pi$ the irreducible representation of $G$
whose highest weight is $-2\pi i\eta$ and by $\pi^*$ its
dual, from Corollary \ref{CorII} it follows that the restriction of $\upsilon$ to
${\cal E}_{\Lambda,I}$ is $\pi^*$. 
  From Theorem  \ref{Thchi1} and Proposition \ref{P1} we deduce

\begin{Thm}\label{BWWey}
Let $\eta$ be an element of ${\frak g}^*$, such that
$2\pi i\eta$ is an 
integral character on ${\frak g}_{\eta}$, and   
$G_{\eta}$ is a maximal torus of $G$. If 
  $\{\psi_t \}$ is the closed Hamiltonian 
isotopy in ${\cal O}_{\eta}$ generated by the 
vector fields $\{X_ {A_t} \}$, then
\begin{equation}\label{Weylf}
\kappa(\psi)=\frac{\chi(\pi^*)(h_1)}{\text{dim}\,\pi},
\end{equation}
 where $h_t\in G$
is the solution to $\Dot h_th_t^{-1}=A_t,\,h_0=e$, and  $\pi$ is the representation of
$G$ whose highest weight is $-2\pi i\eta.$
\end{Thm} 
Now the character $\chi(\pi^*)$ and the dimension $\text{dim}\,\pi$ can be determined by
  Weyl's character formula \cite{JJ00}, and so $\kappa(\psi)$.


\medskip

 {\sc Examples}

{\it 1. Action integral in  flag manifolds.}
Set $D=\text{diag}\,(id_1,\dots,id_n)\in{\frak u}(n)$, with $d_j\in{\Bbb R}$. 
We denote by $p_1<\dots <p_k$ the distinct values of the $d_j$ and by $n_1,\dots,n_k$
the corresponding multiplicities. $D$ determines an element $\eta\in{\frak u}(n)^*$
by the relation $\eta(Y)=\text{tr}(DY)$. The coadjoint orbit ${\cal O}_{\eta}$ is
the flag manifold $U(n)/U(n_1)\times\dots\times U(n_k)$. 
And for $Y=(B_1,\dots,B_k)\in{\frak g}_{\eta}=\bigoplus_j{\frak u}(n_j)$,
$$\eta(Y)=\sum_{j=1}^{k}ip_j\text{tr}(B_j).$$

The manifold ${\cal O}_{\eta}$ depends only on the multiplicities $n_j$. However
 the symplectic form $\omega$, defined in (\ref{omeganu}) depends also on
the $p_j$. The manifold $({\cal O}_{\eta},\omega)$ admits a $U(n)$-invariant 
prequantization if $-2\pi p_j=:m_j\in{\Bbb Z}$, for $j=1,\dots,k.$ In this case the character
$\Lambda$ of $G_{\eta}=\prod_jU(n_j)$ defined by
$$\Lambda(A_1,\dots,A_k)=\prod_{j=1}^k (\text{det}(A_j))^{m_j}$$
has as derivative $2\pi i\eta$. Now   the   symplectic flag manifold 
$({\cal O}_{\eta},\omega)$ is quantizable and it is determined by $(m_1,n_1;\dots;m_k,n_k).$

If $\{g_t\in U(n)\,|\,t\in[0,1] \}$ defines a loop $\varphi$ in 
$\text{Ham}({\cal O}_{\eta})$ by (\ref{g_tx}) , then $g_1gG_{\eta}=gG_{\eta}$, for every $g\in U(n)$.
So $g_1$ is a multiple of the identity; $g_1=zI_n$, with $|z|=1$.
By Theorem \ref{Tmcoador} 
$$\kappa(\varphi)=\Lambda (g_1)=\prod_{j=1}^kz^{m_jn_j}.$$
Thus we have
 
\begin{Prop}
The symplectic flag manifold 
$$({\cal O}_{\eta}=U(n)/U(n_1)\times\dots\times U(n_k),\,\omega)$$
determined by the integers $(m_1,n_1;\dots;m_k,n_k)$
admits a $U(n)$-invariant prequantization. If $\varphi$ is the loop in 
$\text{Ham}({\cal O}_{\eta})$ defined by a family $\{g_t\in U(n)\}$, with  $g_0=I_n$ and 
$g_1=zI_n$, 
then 
$$\kappa(\varphi)=z^a,$$
where $a={\sum_j}m_jn_j.$
\end{Prop}

\medskip   

{\it 2. The invariant $\kappa$ of a Hamiltonian flow in $S^2$.}
 For $G=SU(2)$, if
$$\eta:\pmatrix
ai& w \\
-\bar w& -ai
\endpmatrix\in{\frak su}(2)\mapsto\frac{na}{2\pi}\in{\Bbb R},$$
with $n\in{\Bbb Z}$, then the orbit ${\cal O}_{\eta}=SU(2)/U(1)=S^2$ admits 
and $SU(2)$-invariant quantization and the corresponding character 
$\Lambda$ of $U(1)$ is $\Lambda(z)=z^n$.

Let $E$ be a matrix  of ${\frak su}(2)$, and we assume that $e^E=-\text{Id}$.
If we denote by $\psi_t$ the symplectomorphism of $S^2$ given by
$$\psi_t(q)=\text{exp}(tE)\cdot q,$$
 then the family
$\{\psi_t\}_{t\in[0,1]}$,  is a closed Hamiltonian flow
on the orbit ${\cal O}_{\eta}$. 
 By Theorem \ref{Tmcoador}
$$\kappa(\psi)=\Lambda(e^E)=\Lambda(-\text{Id})=(-1)^n.$$
This result agrees with that one obtained in \cite[Theorem 21]{aV01} by direct calculation.
This value is also obtained in \cite[Example 3.6]{mS00}.

This value can also be deduced from Theorem \ref{BWWey}. 
The Weyl's character formula \cite{JJ00} is very simple for
 the group $SU(2)$; in this case, there is
only one positive root $\alpha$ and
the Weyl group has only two elements. 
We take for $\alpha$ the linear map defined by
$$\alpha(\text{diag}\,(ai,\,-ai))=2ai;$$ 
so $\alpha^{\vee}=\text{diag}\,(1,\,-1)$.
If $n<0$, then $-\lambda:=-2\pi \eta$ is the 
highest weight of a representation $\pi$ of $SU(2)$.
 For $t\in U(1)$, $t^{\lambda}=t^{-n}$ and $t^{\alpha}=t^2$. Therefore (see \cite{JJ00})
$$\text{dim}\, \pi=-n+1 \,\,\,\text{and}\,\,\, \chi_{\pi}(t)=\sum_{k=0}^{-n}t^{-n-2k}.$$
 Hence 
$$\chi_{\pi^*}(h_1)=\chi_{\pi}(-1)=(-n+1)(-1)^n,$$
 and from (\ref{Weylf}) we again obtain  the value $(-1)^n$ for $\kappa(\psi)$.


\end{document}